\documentclass[11pt]{article}

\usepackage{amsmath}
\usepackage{amsfonts}
\usepackage{amssymb}
\usepackage{amscd}
\usepackage{latexsym}
\usepackage{mathrsfs}
\usepackage{amsthm}

\newtheorem{Proposition}{Proposition}

\newtheorem{Remark}{Remark}

\def\C{\mathbb{C}}
\def\R{\mathbb{R}}
\def\d{\mathrm{d}}
\def\F21{{}_2F_1}
\def\e{\varepsilon}
\def\Pv{\mathrm{Pv}}

\begin{document}

\title{Weber-Schafheitlin's type integrals with exponent $1$}

\author{Johannes Kellendonk\,~and\,~Serge Richard}

  \date{\small
    \begin{quote}
      \emph{
    \begin{itemize}
    \item[]
Institut Camille Jordan, CNRS UMR 5208,
Universit\'e Lyon 1, Universit\'e de Lyon,
43 boulevard du 11 novembre 1918, F-69622 Villeurbanne cedex,
France
    \item[]
      \emph{E-mails\:\!:}
      kellendonk@math.univ-lyon1.fr\,~and\,~richard@math.univ-lyon1.fr
    \end{itemize}
      }
    \end{quote}
    March 2008
  }
\maketitle

\begin{abstract}
Explicit formulae for Weber-Schafheitlin's type integrals with exponent $1$ are derived. The results of these integrals are  distributions on $\R_+$.
\end{abstract}

\section{Introduction}
Let $J_\mu$ denote the Bessel function of the first kind and of order $\mu$. An integral of the form
\begin{equation}\label{eqbase}
\int_0^\infty \kappa^\rho \;\!J_\mu (s\kappa) \;\!J_\nu(\kappa)\;\! \d \kappa\ ,
\end{equation}
for suitable $\mu,\nu,\rho  \in \R$ and $s \in \R_+ =(0, \infty)$ is called a Weber-Schafheitlin's integral \cite[Chap.~13.4]{W}. If $\rho$ is strictly less than $1$, the result of this integration is known and can be found in many textbooks or handbooks, see for example \cite{GR,W,Whee}. However, even if the critical case $\rho = 1$ seems quite natural and even necessary for certain calculations on scattering theory \cite{KR}, we have not been able to find it in the literature. Therefore, we provide in this paper the result of \eqref{eqbase} for $\rho = 1$ as well as the result of the related integral
\begin{equation}\label{eqbase2}
\int_0^\infty \kappa \;\!H^{(1)}_\mu (s\kappa) \;\!J_\nu(\kappa)\;\! \d \kappa\ ,
\end{equation}
where $H^{(1)}_\mu$ is the Hankel function of the first kind and of order $\mu$. We emphasize that both results are not functions of the variable $s$ but distributions on $\R_+$. Let us also mention some related recent works \cite{DS,Mil,Mir}. In the first reference, only the special cases $\nu = \pm\mu$ of \eqref{eqbase} and $\nu=\mu$ of \eqref{eqbase2} are explicitly calculated.

\section{The derivation of the integral \eqref{eqbase2}}

Let us start by recalling that for $z \in \C$ satisfying $-\hbox{$\frac{\pi}{2}$}<\arg(z)\leq \pi$ one has \cite[eq.~9.6.4]{AS}~:
\begin{equation*}
H^{(1)}_\mu(z)=\hbox{$\frac{2}{i\pi}$}\;\!
e^{-i\pi\mu/2}\;\!K_\mu(-iz)\ ,
\end{equation*}
where $K_\mu$ is the modified Bessel function of the second kind and of order $\mu$. Moreover, for $\Re(z)>0$ and $\nu + 2 >|\mu|$ the following result holds \cite[Sec.~13.45]{W}~:
\begin{equation*}
\int_0^\infty \kappa \;\!K_\mu (z\kappa) \;\!J_\nu(\kappa)\;\! \d \kappa\ = \hbox{$\frac{\Gamma(\frac{\nu+\mu}{2}+1)\;\!\Gamma (\frac{\nu-\mu}{2}+1)}{\Gamma(\nu+1)}$}\;\! z^{-2-\nu} \;\!\F21\big(\hbox{$\frac{\nu+\mu}{2}+1$},
\hbox{$\frac{\nu-\mu}{2}+1$};\nu+1;-z^{-2}\big) \,
\end{equation*}
where $\F21$ is the Gauss hypergeometric function \cite[Chap.~15]{AS}. Thus, by setting $z=s+i\e$ with $s \in \R_+$ and $\e>0$ one obtains
\begin{eqnarray*}\label{Watson}
& I_{\mu,\nu}(s+i\e):=\int_0^\infty \kappa \;\!H^{(1)}_\mu \big((s+i\e)\kappa\big) \;\!J_\nu(\kappa)\;\! \d \kappa & \\
& = \hbox{$\frac{2}{i\pi}$} \;\!e^{-i\pi\mu/2}\;\!
\hbox{$\frac{\Gamma(\frac{\nu+\mu}{2}+1)\;\!\Gamma
(\frac{\nu-\mu}{2}+1)}{\Gamma(\nu+1)}$}\;\!
(-is +\e)^{-2-\nu} \;\!\F21\big(\hbox{$\frac{\nu+\mu}{2}$}+1,
\hbox{$\frac{\nu-\mu}{2}$}+1;\nu+1;(s+i\e)^{-2}\big)\ .&
\end{eqnarray*}
Taking into account Equality 15.3.3 of \cite{AS} one can isolate from the $\F21$-function a factor which is singular for $s=1$ when $\e$ goes to $0$:
\begin{equation*}
\F21\big(\hbox{$\frac{\nu+\mu}{2}$}+1,\hbox{$\frac{\nu-\mu}{2}$}+1;
\nu+1;(s+i\e)^{-2}\big) = \hbox{$\frac{1}{1-(s+i\e)^{-2}}$}\;
\F21\big(\hbox{$\frac{\nu+\mu}{2}$},\hbox{$\frac{\nu-\mu}{2}$};
\nu+1;(s+i\e)^{-2}\big)\ .
\end{equation*}
Furthermore, by inserting the equalities
\begin{equation*}
\hbox{$\frac{1}{1-(s+i\e)^{-2}}$} =
-(s+i\e)^2\;\!\hbox{$\frac{1}{s}$}
\;\!\hbox{$\frac{1}{\big(\frac{(1+\e^2)}{s}-s\big)-
2i\e}$}
\end{equation*}
and
\begin{equation*}
(-is + \e)^{-2-\nu}=
-e^{i\pi\nu/2}(s+i\e)^{-2-\nu}\ ,
\end{equation*}
in these expressions, one finally obtains that $I_{\mu,\nu}(s+i\e)$ is equal to
\begin{equation}\label{Tokyo1}
\hbox{$\frac{2}{i\pi}$} \;\!e^{i\pi(\nu-\mu)/2}\;\!\hbox{$\frac{1}{s}$}\;\!
\hbox{$\frac{(s+i\e)^{-\nu}}{\big(\frac{(1+\e^2)}{s}-s\big)-
2i\e}$}\;\! \hbox{$\frac{\Gamma(\frac{\nu+\mu}{2}+1)\;\!\Gamma
(\frac{\nu-\mu}{2}+1)}{\Gamma(\nu+1)}$}\;\!
\F21\big(\hbox{$\frac{\nu+\mu}{2}$},\hbox{$\frac{\nu-\mu}{2}$};
\nu+1;(s+i\e)^{-2}\big) \ .
\end{equation}

We are now ready to study the $\e$-behaviour of each of the above
terms. For the particular choice of the three parameters  $\hbox{$\frac{\nu+\mu}{2}$}$, $\hbox{$\frac{\nu-\mu}{2}$}$ and $\nu+1$, the map
\begin{equation*}
z \mapsto
\F21\big(\hbox{$\frac{\nu+\mu}{2}$},\hbox{$\frac{\nu-\mu}{2}$};
\nu+1;z\big)\ ,
\end{equation*}
which is holomorphic in the cut complex plane $\C\setminus [1,\infty)$, extends continuously to  $[1,\infty)$. The limits from above and below yield generally two different continuous functions, and by convention the hypergeometric function on $[1,\infty)$ is the limit obtained from below. Since $\Im\big((s+i\e)^{-2}\big)<0$, the $\F21$-factor in \eqref{Tokyo1} converges to $ \F21\big(\hbox{$\frac{\nu+\mu}{2}$}, \hbox{$\frac{\nu-\mu}{2}$};\nu+1;s^{-2}\big)$ as $\e \to 0$, uniformly in $s$ on any compact subset of $\R_+$.

For the other factors, let us observe that $(s+i\e)^{-\nu}$ converges to $s^{-\nu}$ as $\e\to 0$ uniformly in $s$ on any compact subset of $\R_+$. Furthermore, it is known that
\begin{equation}\label{Tokyo2}
\lim_{\e \to 0} \hbox{$\frac{1}{\big(\frac{(1+\e^2)}{s}-s\big)-
2i\e}$} = \Pv\;\!\big(\hbox{$\frac{1}{\frac{1}{s}-s}$}\big)
+i\hbox{$\frac{\pi}{2}$}\delta(s-1)\ ,
\end{equation}
where the convergence has to be understood in the sense of distributions on $\R_+$. In the last expression, $\delta$ is the Dirac measure centered at $0$ and $\Pv$ denotes the principal value integral. By collecting all these results one can prove:

\begin{Proposition}\label{surImunu}
For any $\mu,\nu\in \R$ satisfying $\nu + 2 >|\mu|$ and $s \in \R_+$ one has
\begin{eqnarray}\label{horreur1}
& \int_0^\infty \kappa \;\!H^{(1)}_\mu (s\kappa) \;\!J_\nu(\kappa)\;\! \d \kappa\ = e^{i\pi(\nu-\mu)/2}\;\!\delta(s-1)  &\\
\nonumber & + \hbox{$\frac{2}{i\pi}$} \;\!e^{i\pi(\nu-\mu)/2}\;\!\Pv \big(\hbox{$\frac{1}{\frac{1}{s}-s}$}\big)\;\!
\;\!\hbox{$\frac{s^{-\nu}}{s}$}\hbox{$\frac{\Gamma(\frac{\nu+\mu}{2}+1)\;\!\Gamma(\frac{\nu-\mu}{2}
+1)}{\Gamma(\nu+1)}$}\;\!
\F21\big(\hbox{$\frac{\nu+\mu}{2}$},\hbox{$\frac{\nu-\mu}{2}$}
;\nu+1;s^{-2}\big)&
\end{eqnarray}
as an equality between two distributions on $\R_+$.
\end{Proposition}

\begin{Remark}\label{rem2}
A priori, the second term in the r.h.s.~is not well defined, since it is the product of the distribution $\Pv \big(\hbox{$\frac{1}{\frac{1}{s}-s}$}\big)$ with a function which is not smooth, or at least differentiable at $s=1$. However, by using the development of the hypergeometric function in a neighbourhood of $1$ \cite[eq.~15.3.11]{AS} it is easily observed that
\begin{equation}\label{moinsenerve}
\hbox{$\frac{\Gamma(\frac{\nu+\mu}{2}+1)\;\!\Gamma(\frac{\nu-\mu}{2}
+1)}{\Gamma(\nu+1)}$}\;\!
\F21\big(\hbox{$\frac{\nu+\mu}{2}$},\hbox{$\frac{\nu-\mu}{2}$}
;\nu+1;s^{-2}\big) = 1 +(s-1)\;\!h(s)
\end{equation}
with a function $h$ that belongs to $L^1_{\rm loc}(\R_+)$. Thus, for any $\alpha \in \R$, the second term in the r.h.s.~of \eqref{horreur1} is equal to
\begin{equation}\label{correction}
\hbox{$\frac{2}{i\pi}$} \;\!e^{i\pi(\nu-\mu)/2}\;\!\Big[
s^{\alpha}\;\!\Pv \big(\hbox{$\frac{1}{\frac{1}{s}-s}$}\big) +
\hbox{$\frac{s^{\alpha}}{\frac{1}{s}-s}$}
\Big(\hbox{$\frac{s^{-\nu-\alpha}}{s}$}\hbox{$\frac{\Gamma(\frac{\nu+\mu}{2}+1)
\;\!\Gamma(\frac{\nu-\mu}{2}+1)}{\Gamma(\nu+1)}$}\;\!
\F21\big(\hbox{$\frac{\nu+\mu}{2}$},\hbox{$\frac{\nu-\mu}{2}$}
;\nu+1;s^{-2}\big) -1\Big)\Big] \ ,
\end{equation}
with the second term in $L^1_{\rm loc}(\R_+)$. Clearly, this distribution is now well defined. The parameter $\alpha$ has been added because it may be useful in certain applications.
\end{Remark}

\begin{Remark}
To describe the singularity at $s=1$, the decomposition \eqref{correction} of the second term of the r.h.s.~of \eqref{horreur1} is certainly valuable. However, it seems to us that this decomposition is less useful if one needs to control the behaviour at $s=0$.
\end{Remark}

\begin{Remark}
In the special case $\mu=\pm \nu$, the hypergeometric function is equal to $1$, and thus the r.h.s.~simplifies drastically.
\end{Remark}

\begin{proof}
a) For any $\e>0$, let us define the function
\begin{equation*}
p_\e : \R_+ \ni s \mapsto \hbox{$\frac{1}{\big(\frac{(1+\e^2)}{s}-s\big)-
2i\e}$}\in \C\ ,
\end{equation*}
and the function $\R_+\ni s \to q_\e(s) \in \C$ by
\begin{equation*}
q_\e(s):=\hbox{$\frac{2}{i\pi}$}
\;\!e^{i\pi(\nu-\mu)/2}\;\!\hbox{$\frac{(s+i\e)^{-\nu}}{s}$}
\hbox{$\frac{\Gamma(\frac{\nu+\mu}{2}+1)\;\!\Gamma(\frac{\nu-\mu}{2}
+1)}{\Gamma(\nu+1)}$}\;\!
\F21\big(\hbox{$\frac{\nu+\mu}{2}$},\hbox{$\frac{\nu-\mu}{2}$}
;\nu+1;(s+i\e)^{-2}\big)\ .
\end{equation*}
Clearly, one gets from \eqref{Tokyo1} that $p_\e(s) \;\! q_\e(s)=I_{\mu,\nu}(s+i\e)$ and from the above remarks that
\begin{equation*}
\lim_{\e\to 0}q_\e(s) =\hbox{$\frac{2}{i\pi}$} \;\!
e^{i\pi(\nu-\mu)/2}\hbox{$\frac{s^{-\nu}}{s}$}\;\!\hbox{$\frac{\Gamma(\frac{\nu+\mu}{2}+1)\;\!
\Gamma(\frac{\nu-\mu}{2}+1)}{\Gamma(\nu+1)}$}
\F21\big(\hbox{$\frac{\nu+\mu}{2}$},\hbox{$\frac{\nu-\mu}{2}$}
;\nu+1;s^{-2}\big)=:q_0(s)\ ,
\end{equation*}
the convergence being uniform in $s$ on any compact subset of
$\R_+$. Furthermore, it follows from \cite[eq.~15.1.20]{AS} that $q_0(1)=\hbox{$\frac{2}{i\pi}$} \;\!e^{i\pi(\nu-\mu)/2}$.

b) Let $g \in C^\infty_c(\R_+)$, {\it i.e.}~$g$ is a smooth function with compact support in $\R_+$, and set
\begin{equation*}
M_g:=\sup_{\e \in [0,1]}\int_{\R_+} |q_\e(s)\;\!g(s)|
\;\!\d s\ .
\end{equation*}
For any $\eta>0$, one can then choose a compact subset $K_\eta$ of $[0,\infty)$ such that
\begin{equation*}
\sup_{\e \in [0,1]}\ \sup_{s \in \R_+\setminus
K_\eta}|p_\e(s)|\leq \hbox{$\frac{\eta}{6 M_g}$}\ ,
\end{equation*}
which implies that $\int_{\R_+\setminus K_\eta}\big|q_\e(s)\;\!
p_{\e'}(s)\;\!g(s)\big|\d s\leq \hbox{$\frac{\eta}{6}$}$ for
all $\e,\e' \in [0,1]$.

c) Finally, one has
\begin{eqnarray*}
\int_{\R_+}q_\e(s)\;\! p_\e(s)\;\!g(s)\;\! \d s & = &
\int_{\R_+}q_0(s)\;\!p_\e(s)\;\!g(s)\;\! \d s +
\int_{\R_+}\big(q_\e(s)-q_0(s)\big)\;\! p_\e(s)\;\!g(s)\;\!
\d s \\
& = & \int_{\R_+}p_\e(s)\;\!q_0(s)\;\!g(s)\;\! \d s +
\int_{K_\eta}\big(q_\e(s)-q_0(s)\big)\;\!p_\e(s)\;\!g(s)\;\!
\d s \\
&& + \int_{\R_+\setminus
K_\eta}\big(q_\e(s)-q_0(s)\big)\;\!p_\e(s) \;\!g(s)\;\! \d
s \ .
\end{eqnarray*}
For $\e\to 0$, the first term on the r.h.s.~converges to
\begin{eqnarray*}
e^{i\pi(\nu-\mu)/2} g(1) +
\hbox{$\frac{2}{i\pi}$} \;\!e^{i\pi(\nu-\mu)/2}
\int_{\R_+} \Pv\big(\hbox{$\frac{1}{\frac{1}{s}-s}$}\big)
\hbox{$\frac{s^{-\nu}}{s}$}\hbox{$\frac{\Gamma(\frac{\nu+\mu}{2}+1)\;\!
\Gamma(\frac{\nu-\mu}{2}+1)}{\Gamma(\nu+1)}$}\;\!
\F21\big(\hbox{$\frac{\nu+\mu}{2}$},\hbox{$\frac{\nu-\mu}{2}$}
;\nu+1;s^{-2}\big)\;\!g(s)\d s.
\end{eqnarray*}
Indeed, the convergence of \eqref{Tokyo2} holds in the sense of dis\-tri\-butions not only on smooth functions on $\R_+$ with compact support, but also on the product $q_0 \;\!g$ of the non-smooth function $q_0$ with the smooth function $g$. This can easily be obtained by using the development given in \eqref{moinsenerve}.
Furthermore, one has
\begin{equation*}
\big|\int_{K_\eta}\big(q_\e(s)-q_0(s)\big)\;\!p_\e(s)
\;\!g(s)\;\! \d s \big| \leq \sup_{s \in
K_\eta\cap \;\mathrm{ supp }\; g}|q_\e(s)-q_0(s)|\int_{K_\eta} |p_\e(s)
\;\!g(s)| \d s
\end{equation*}
which is less that $\hbox{$\frac{\eta}{3}$}$ for $\e$ small
enough since $q_\e-q_0$ converges uniformly to $0$ on any com\-pact
subset of $\R_+$. And finally, from the choice of $K_\eta$ one has
\begin{equation*}
\big| \int_{\R_+\setminus
K_\eta}\big(q_\e(s)-q_0(s)\big)\;\!p_\e(s) \;\!g(s)\;\! \d
s \big| \leq \hbox{$\frac{\eta}{3}$}\ .
\end{equation*}
Since $\eta$ is arbitrary, one has thus obtained that the map $s \mapsto p_\e(s)\;\!q_\e(s)$ converges in the sense of distributions on $\R_+$ to the distribution given by the r.h.s.~term of the statement of the proposition.

\end{proof}

\begin{Remark}
In the previous proof, the Lebesgue measure on $\R_+$ has been used for the evaluation of the distribution on a smooth function with compact support in $\R_+$. Let us notice that the same result holds if the measure $\hbox{$\frac{\d s}{s}$}$ is chosen instead of the Lebesgue measure.
\end{Remark}

\section{The derivation of the integral \eqref{eqbase}}

In the next proposition, the function $H^{(1)}_\mu$ of the previous statement is replaced by the Bessel function $J_\mu$. Since $H_\mu^{(1)} = J_\mu + iY_\mu$ with $J_\mu$ and $Y_\mu$ real on $\R_+$, taking the real part of both side of \eqref{horreur1} would lead to the result. However, since the real and the imaginary parts of the Gauss hypergeometric function are not very explicit, we prefer to sketch an independent proof.

\begin{Proposition}\label{JmuJnu}
For any $\mu,\nu \in \R$ satisfying $\nu+2>|\mu|$ and $\mu+2>|\nu|$, and $s \in \R_+$
one has
\begin{eqnarray*}
&\int_0^\infty \kappa \;\!J_\mu (s\kappa) \;\!J_\nu(\kappa)\;\! \d \kappa\ = \cos(\pi(\nu -\mu)/2)\;\!\delta(s-1) & \\
&+
\hbox{$\frac{2}{\pi}$} \;\!\sin(\pi(\nu -\mu)/2) \;\!
\Pv\big(\hbox{$\frac{1}{\frac{1}{s}-s}$}\big)\cdot
\left\{
\begin{array}{ll}
\hbox{$\frac{s^\mu}{s}$}\hbox{$\frac{\Gamma(\frac{\mu+\nu}{2}+1)\;\!
\Gamma(\frac{\mu-\nu}{2}+1)}{\Gamma(\mu+1)}$}\;\!
\F21\big(\hbox{$\frac{\mu+\nu}{2}$},\hbox{$\frac{\mu-\nu}{2}$};\mu+1;s^2\big)
&\mbox{ if } s\leq 1,\\
\hbox{$\frac{s^{-\nu}}{s}$}\hbox{$\frac{\Gamma(\frac{\nu+\mu}{2}+1)\;\!
\Gamma(\frac{\nu-\mu}{2}+1)}{\Gamma(\nu+1)}$}\;\!
\F21\big(\hbox{$\frac{\nu+\mu}{2}$},\hbox{$\frac{\nu-\mu}{2}$};\nu+1;s^{-2}\big)
&\mbox{ if }s> 1,
\end{array}
\right.&
\end{eqnarray*}
as an equality between two distributions on $\R_+$.
\end{Proposition}

\begin{Remark}
We refer to Remark \ref{rem2} for a discussion on the
fact that the distribution corresponding to the second term on the r.h.s.~is well defined.
\end{Remark}

\begin{proof}
a) Even if the proof is very similar to the previous one, two additional observations have to be taken into account: 1) The map $z\mapsto \F21(\alpha,\beta;\gamma;z)$ is real when $z$ is restricted to the interval $[0,1)$, 2) There exists a simple relation between $\Re(I_{\mu,\nu})$ and $\Re(I_{\nu,\mu})$. More precisely, for any $s \in \R_+\setminus\{1\}$ one has
\begin{eqnarray*}
\Re\big(I_{\mu,\nu}(s)\big) &=&\int_0^\infty \kappa \;\!J_\mu (s\kappa) \;\!J_\nu(\kappa)\;\! \d \kappa \\
&=& s^{-2} \int_0^\infty \kappa \;\!J_\mu (\kappa) \;\!J_\nu(s^{-1}\kappa)\;\! \d \kappa \\
&=& s^{-2}\;\!\Re\big(I_{\nu,\mu}(s^{-1})\big)\ .
\end{eqnarray*}
Thus, the main trick of the proof is to use the previous expressions for $s\in (1,\infty)$ since the contribution of the $\F21$-function in \eqref{horreur1} is then real, and to obtain similar formulae below for $s\in(0,1]$. However, some care has to be taken because of the Dirac measure at $1$ and of the principal value integral also centered at $1$.

So, for any $s\in (0,1]$ and $\e>0$, let us set $z:=s-i\e$. Since the conditions $-\hbox{$\frac{\pi}{2}$}<\arg(z^{-1})\leq \pi$, $\Re(-iz^{-1})>0$ and $\mu+2>|\nu|$ hold, one can obtain the analog of \eqref{Watson}:
\begin{equation*}
\hbox{$\frac{1}{z^2}$}\int_0^\infty \kappa \;\!H^{(1)}_\nu (z^{-1}\kappa) \;\!J_\mu(\kappa)\;\! \d \kappa  =
-\hbox{$\frac{2}{i\pi}$}\;\!e^{-i\pi(\nu-\mu)/2}\;\!
\hbox{$\frac{1}{\big(\frac{(1+\e^2)}{s}-s\big)+
2i\e}$}\;\!m_\e(s) \ ,
\end{equation*}
with $m_\e:(0,1]\to\C$ the function defined by
\begin{equation}\label{defdem}
m_\e(s):=\hbox{$\frac{(s-i\e)^\mu}{s}$}\;\!
\hbox{$\frac{\Gamma(\frac{\mu+\nu}{2}+1)\;\!
\Gamma (\frac{\mu-\nu}{2}+1)}{\Gamma(\mu+1)}$}\;\! \F21\big(\hbox{$\frac{\mu+\nu}{2}$},
\hbox{$\frac{\mu-\nu}{2}$};\mu+1;(s-i\e)^2\big)\ .
\end{equation}

b) We now collect both functions, for $s\leq 1$ and $s>1$. For that purpose, let $l_\e:\R_+ \to \C$ be the function defined for $s \in (0,1]$ by
\begin{equation*}
l_\e(s)= -\hbox{$\frac{2}{i\pi}$}\;\!e^{-i\pi(\nu-\mu)/2}\;\!
\hbox{$\frac{1}{\big(
\frac{(1+\e^2)}{s}-s\big)+ 2i\e}$}
\end{equation*}
and for $s \in (1,\infty)$ by
\begin{equation*}
l_\e(s)= \hbox{$\frac{2}{i\pi}$}\;\!e^{i\pi(\nu-\mu)/2}\;\!
\hbox{$\frac{1}{\big(\frac{(1+\e^2)}{s}-
s\big)- 2i\e}$}\ .
\end{equation*}
By selecting only the real part of $l_\e$, one obtains that $\Re(l_\e)$ converges in the sense of distributions on $\R_+$ as $\e\to 0$ to the distribution:
\begin{equation*}
\cos(\pi(\nu -\mu)/2)\;\!\delta(s-1)
 + \hbox{$\frac{2}{\pi}$} \;\!\sin(\pi(\nu -\mu)/2) \;\! \;\!\Pv\;\!\big(\hbox{$\frac{1}{\frac{1}{s}-s}$}\big)\ .
\end{equation*}
Furthermore, let $m_\e:\R_+\to\C$ be the function defined for $s \in (0,1]$ by
\eqref{defdem} and for $s\in (1,\infty)$ by
\begin{equation*}
\hbox{$\frac{(s+i\e)^{-\nu}}{s}$}
\hbox{$\frac{\Gamma(\frac{\nu+\mu}{2}+1)\;\!\Gamma(\frac{\nu-\mu}{2}
+1)}{\Gamma(\nu+1)}$}\;\!
\F21\big(\hbox{$\frac{\nu+\mu}{2}$},\hbox{$\frac{\nu-\mu}{2}$}
;\nu+1;(s+i\e)^{-2}\big).
\end{equation*}
As $\e$ goes to $0$, this function converges uniformly on any compact subset of $\R_+$ to a continuous real function $m_0$. Furthermore, this function takes the value $1$ for $s=1$. Indeed, these properties of $m_0$ follow from the facts that the $\F21$-functions are real on the $[0,1]$, that the convergences from above or below this interval give the same values, and that the normalization factors have been suitably chosen. The function $m_0$ is the one given after the curly bracket in the statement of the proposition.

c) The remaining part of the proof can now be mimicked from part c) of the proof of the previous proposition. Therefore, we simply refer to this paragraph and omit the end of the present proof.
\end{proof}

\end{document}